\documentclass[12pt,twoside]{article}

\usepackage{amssymb,amsmath}

\usepackage[active]{srcltx} 

\newtheorem{theorem}{Theorem}[section]

\numberwithin{equation}{section}

\newcommand{\N}{\mathbb{N}}

\newcommand{\C}{\mathbb{C}}
\newcommand{\dis}{\displaystyle}

\newcommand{\db}{\rule[.05in]{.09in}{.10in}} 

\begin{document}

\pagestyle{myheadings}

\markboth{Igor E. Pritsker}{Approximation of conformal mapping via
the Szeg\H{o} kernel method}%

\title{Approximation of conformal mapping via the Szeg\H{o} kernel method}%
\author{Igor E. Pritsker\thanks{This material is based upon work supported
by the National Science Foundation under Grant No. 9996410, and by
the National Security Agency under Grant No. MDA904-03-1-0081.}}%

\date{}%

\maketitle

\centerline{\it Dedicated to the memory of Professor D. Gaier}

\begin{abstract}

We study the uniform approximation of the canonical conformal
mapping, for a Jordan domain onto the unit disk, by polynomials
generated from the partial sums of the Szeg\H{o} kernel expansion.
These polynomials converge to the conformal mapping uniformly on
the closure of any Smirnov domain. We prove estimates for the rate
of such convergence on domains with piecewise analytic boundaries,
expressed through the smallest exterior angle at the boundary.
Furthermore, we show that the rate of approximation on compact
subsets inside the domain is essentially the square of that on the
closure. Two standard applications to the rate of decay for the
contour orthogonal polynomials inside the domain, and to the rate
of locally uniform convergence of Fourier series are also given.
\\
\\
{\bf Keywords.} Conformal mapping, Szeg\H{o} kernel, orthogonal
polynomials, Fourier series.
\\
\\
{\bf 2000 MCS.} Primary 30C40, 30E10; Secondary 41A10, 30C30.

\end{abstract}


\section{Convergence of the Szeg\H{o} kernel expansion and
approximation of conformal maps}

Let $G$ be a Jordan domain in the complex plane. There are two
well known kernel methods used for approximation of the canonical
conformal mappings of $G$ onto a disk. The Bergman kernel method
is associated with the $L_2$ spaces and orthogonal polynomials
with respect to the area measure, while the Szeg\H{o} kernel
method is based on the inner product and orthogonal polynomials
with respect to the arclength measure on the boundary of $G$ (see
Gaier \cite{Ga1}, Smirnov and Lebedev \cite{SL}). The
approximations related to the area orthogonality approach were
first introduced by Bieberbach \cite{Bi} through an extremal
problem for polynomials. The first result about the uniform
convergence of the Bieberbach polynomials was proved by Keldysh
\cite{Ke2}, for domains with sufficiently smooth boundaries. The
uniform convergence of the Bieberbach polynomials has been
extensively studied since then (see Gaier \cite{Ga2} and
references therein). A selection of further results on this
subject is in the papers of Andrievskii \cite{An}, Gaier
\cite{Ga3}-\cite{Ga5}, Andrievskii and Gaier \cite{AG}, and
Andrievskii and Pritsker \cite{AP}. In contrast, the method based
on the Szeg\H{o} kernel did not receive such a comprehensive
attention. The goal of this paper is to show that the Szeg\H{o}
kernel method is also very useful for the uniform approximation of
conformal mappings.

Suppose that $G$ has rectifiable boundary $L$ of length $l.$ We
consider the Smirnov spaces $E_p(G), \ 1\le p<\infty,$ of analytic
functions in $G$, whose boundary values satisfy
\[
\|f\|_p = \left( \frac{1}{l}\int_{L} |f(z)|^p | dz |
\right)^{1/p}< \infty
\]
(see Duren \cite{Du}, Smirnov and Lebedev \cite{SL}). Our interest
is focused on the Hilbert space $E_2(G)$, equipped with the inner
product
\[
(f,g) := \frac{1}{l} \int_{L} f(z) \overline{g(z)} |dz|, \quad f,g
\in E_2(G).
\]
Polynomials are dense in $E_p(G), \ 1\le p<\infty,$ if and only if
$G$ is a Smirnov domain \cite{Du}. Keldysh and Lavrentiev
\cite{KL} characterized Smirnov domains by the property that, for
a conformal mapping $\psi$ of the unit disk $D$ onto $G$,
$\log|\psi|$ is represented by the Poisson integral of its
boundary values. Although no complete geometric description of
Smirnov domains is known, this class is sufficiently wide. In
particular, it contains (in the decreasing order of generality)
all  Ahlfors-regular domains, Lavrentiev (cord-arc) domains,
Lipschitz domains, domains with bounded boundary rotation (Radon
domains), piecewise smooth and smooth domains (cf. Pommerenke
\cite[Chap. 7]{Po}). Applying Gram-Schmidt orthonormalization to
monomials $\{ z^n \}_{n =0}^{\infty}$ in the Smirnov domain $G$,
we obtain a complete orthonormal system of polynomials $\{ p_n(z)
\}_{n =0}^{\infty}$ in $E_2(G)$ (see Szeg\H{o} \cite{Sz}). Next,
we introduce the Szeg\H{o} kernel
\begin{equation} \label{1.1}
K(z,\zeta)=\sum_{k=0}^{\infty} \overline{p_k(\zeta)} p_k(z), \quad
z, \zeta \in G,
\end{equation}
where convergence of this bilinear series is uniform in $z$ and
$\zeta$ on compact subsets in $G$ \cite{Sz}, \cite{SL},
\cite{Ga1}. The importance of the Szeg\H{o} kernel lies in its
reproducing property
\begin{equation} \label{1.2}
f(\zeta)= \frac{1}{l} \int_L f(z)\overline{K(z,\zeta)} |dz|, \quad
\zeta \in G,
\end{equation}
which holds for any $f \in E_2(G).$ Equivalently, every $f \in
E_2(G)$ can be represented by its Fourier series
\begin{equation} \label{1.3}
f(\zeta)= \sum_{k=0}^{\infty} a_k p_k(\zeta) = \sum_{k=0}^{\infty}
\left(\frac{1}{l} \int_L f(z)\overline{p_k(z)} |dz|\right)
p_k(\zeta),
\end{equation}
convergent in $E_2(G)$ norm and, consequently, locally uniformly
convergent in $G$ (see Chapter 4 of \cite{SL}).

It is well known that the Szeg\H{o} kernel is closely connected
with the canonical conformal mapping $\varphi$ of $G$ onto the
unit disk $D$:
\begin{equation} \label{1.4}
K(z,\zeta)=\frac{l}{2\pi} \sqrt{\varphi'(z)\varphi'(\zeta)},
\end{equation}
where $\varphi(\zeta)=0$ and $\varphi'(\zeta)>0$ (cf. \cite{SL}
and \cite{Sz}). Hence we have that
\begin{equation} \label{1.5}
\varphi'(\zeta)=\frac{2\pi}{l}\ K(\zeta,\zeta) =\frac{2\pi}{l}
\sum_{k=0}^{\infty} |p_k(\zeta)|^2,
\end{equation}
by \eqref{1.1} and \eqref{1.4}. It follows that
\begin{equation} \label{1.6}
\varphi'(z)=\frac{2\pi}{l}\, \frac{(K(z,\zeta))^2}
{K(\zeta,\zeta)} = \frac{2\pi}{l}\
\frac{\left(\dis\sum_{k=0}^{\infty} \overline{p_k(\zeta)}
p_k(z)\right)^2}{\dis\sum_{k=0}^{\infty} |p_k(\zeta)|^2}, \quad z
\in G,
\end{equation}
where $\zeta\in G$ is regarded as a fixed point. We now introduce
the following sequence of approximating polynomials:
\begin{equation} \label{1.7}
J_{2n+1}(z)=\frac{2\pi}{l}\ \frac{\dis\int_{\zeta}^z
\left(\sum_{k=0}^n \overline{p_k(\zeta)} p_k(t)\right)^2
dt}{\dis\sum_{k=0}^n |p_k(\zeta)|^2}, \quad n \in \N.
\end{equation}
Note that the degree of $J_{2n+1}(z)$ is $2n+1.$ The sequence
$\{J_{2n+1}\}_{n=0}^{\infty}$ converges to $\varphi$ uniformly on
compact subsets of $G$, which is inherited from the partial sums
of the Szeg\H{o} kernel. Similar approximating polynomials, but
with a different normalization, were introduced via an extremal
problem for any $E_p(\Gamma),\ 1 \le p < \infty,$ by Keldysh and
Lavrentiev (see \cite{Ke1} and \cite{KL}). They developed the
ideas of Julia \cite{Ju}, who earlier considered the same extremal
problem for the conformal mapping. Further study of the
convergence properties in $E_2(G)$ is due to Warschawski \cite{Wa}
(also see Gaier \cite{Ga1} for a survey). Convergence questions
for general Fourier expansions in contour orthogonal polynomials
were considered in Rosenbloom and Warschawski \cite{RW}, Smirnov
and Lebedev \cite{SL} and Suetin \cite{Su}. However, all of these
studies impose quite strict smoothness assumptions on the boundary
of $G$. We prove the first uniform convergence results for domains
with corners, and give explicit rates of approximation in terms of
the geometric properties of $G$.

We start with a general estimate for the uniform (sup) norm of the
approximation error on $\overline{G}.$ A similar result was proved
by Warschawski \cite{Wa}, but in a somewhat different form (see
also \cite[pp. 130-131]{Ga1}).

\begin{theorem} \label{thm1.1}
Let $G$ be a Smirnov domain, with $\zeta\in G$ fixed. If
$\varphi:G\to D$ is a conformal mapping normalized by
$\varphi(\zeta)=0$ and $\varphi'(\zeta)>0,$ then
\begin{equation} \label{1.8}
\|\varphi-J_{2n+1}\|_{\infty} \le 8\pi
\left\|K(\cdot,\zeta)-\sum_{k=0}^n \overline{p_k(\zeta)}
p_k(\cdot)\right\|_2 \to 0 \quad\mbox{as } n\to\infty.
\end{equation}
\end{theorem}

Recall that $\sum_{k=0}^n \overline{p_k(\zeta)} p_k(z)$ is the
Fourier sum for $K(z,\zeta)$, i.e., it is the best $E_2(G)$
approximation from the subspace of polynomials of degree $n.$
Hence we can give an upper estimate for the rate of convergence by
appropriately choosing a sequence of approximating polynomials for
$K(z,\zeta)$. The rate of convergence necessarily depends on the
geometric properties of the domain $G$. We consider a class of
domains with piecewise analytic boundaries, which is important in
applications. An analytic arc is defined as the image of a segment
under a mapping that is conformal in an open neighborhood of the
segment. Thus a domain has piecewise analytic boundary if it is
bounded by a Jordan curve consisting of a finite number of
analytic arcs.

\begin{theorem} \label{thm1.2}
Let $\partial G$ be piecewise analytic, with the smallest exterior
angle $\lambda\pi,\ 0<\lambda<2,$ at the junction points of the
analytic arcs. If $\zeta\in F$, where $F\subset G$ is compact,
then
\begin{equation} \label{1.9}
\left\|K(\cdot,\zeta)-\sum_{k=0}^n \overline{p_k(\zeta)}
p_k(\cdot)\right\|_2 \le C_1(G,F)\
n^{-\frac{\lambda}{4-2\lambda}}, \quad n\in\N.
\end{equation}
Here, the constant $C_1(G,F)>0$ depends only on $G$ and $F$.
\end{theorem}

Combining Theorems \ref{thm1.1} and \ref{thm1.2}, we obtain the
main result on the uniform approximation of conformal mappings.

\begin{theorem} \label{thm1.3}
Let $\partial G$ be piecewise analytic, with the smallest exterior
angle $\lambda\pi,\ 0<\lambda<2,$ at the junction points of the
analytic arcs. Suppose that $\zeta\in F$, where $F\subset G$ is
compact. For the conformal mapping $\varphi:G\to D$, normalized by
$\varphi(\zeta)=0$ and $\varphi'(\zeta)>0,$ we have
\begin{equation} \label{1.10}
\|\varphi-J_{2n+1}\|_{\infty} \le C_2(G,F)\
n^{-\frac{\lambda}{4-2\lambda}}, \quad n\in\N,
\end{equation}
where the constant $C_2(G,F)>0$ depends only on $G$ and $F$.
\end{theorem}

It is interesting to compare the convergence properties of our
sequence $\{J_{2n+1}\}_{n=0}^{\infty}$ with those of the
Bieberbach polynomials $\{B_n\}_{n=0}^{\infty}$. Gaier \cite{Ga4}
proved for domains with piecewise analytic boundaries that
\[
\|\varphi-B_n\|_{\infty} \le C_3\ \log{n}\
n^{-\frac{\lambda}{2-\lambda}}, \quad n\ge 2.
\]
Later, Andrievskii and Gaier \cite{AG} replaced $\log{n}$ by
$\sqrt{\log{n}},$ and relaxed the imposed geometric condition to
piecewise quasianalytic boundary. Their estimate for the rate of
uniform convergence of the Bieberbach polynomials remains the most
precise known. Although \eqref{1.10} gives a slower rate of
convergence, the polynomials $J_{2n+1}$ have some advantages over
the Bieberbach polynomials. They are free from the convergence
anomalies exhibited by Keldysh's example \cite{Ke2}, where only
one singular point, at otherwise very smooth boundary, destroys
the uniform convergence of the Bieberbach polynomials. (One can
find further information about this example of Keldysh in
\cite{AP}.) In fact, Smirnov domains include any imaginable domain
arising in numerical applications, guaranteeing the uniform
convergence of $\{J_{2n+1}\}_{n=0}^{\infty}$ by Theorem
\ref{thm1.1}. In addition, this sequence is easier to generate
numerically, because the inner products defined by the boundary
integrals are easier to compute than the area inner products in
Gram-Schmidt orthonormalization.

We show that the rate of convergence for $J_{2n+1}$ on compact
subsets of $G$ is better than on the whole domain, i.e., it is
essentially squared comparing to \eqref{1.10}.

\begin{theorem} \label{thm1.4}
If the conditions of Theorem \ref{thm1.3} are satisfied, then
\begin{equation} \label{1.11}
\max_{z\in F} |\varphi(z)-J_{2n+1}(z)| \le C_4(G,F)\
n^{-\frac{\lambda}{2-\lambda}}, \quad n\in\N,
\end{equation}
where $C_4(G,F)>0$ depends only on $G$ and $F$.
\end{theorem}

Gaier \cite{Ga5} posed the question of possible improvement in
locally uniform convergence rates for the Bieberbach polynomials
and other approximations for conformal maps. The above theorem
provides a partial answer for his question in the case of the
Szeg\H{o} kernel method. Furthermore, it is possible to give
similar improvements for the Bieberbach polynomials too, by
following the ideas of this paper.

Clearly, it was not our goal to achieve the highest possible level
of generality here. Thus Theorem \ref{thm1.2} is true for domains
with piecewise quasianalytic boundaries \cite{AG} (therefore, all
other results are valid for these domains too). One only needs to
fill in a number of technical details on the behavior of
$\sqrt{\varphi'}$ near corners, in our proof, to reach this
conclusion. We are also able to prove Theorems
\ref{thm1.2}-\ref{thm1.4} for Lavrentiev (cord-arc) domains, with
the rates of convergence of the order $n^{-\gamma},\ \gamma>0.$
Finally, we have analogues of these results in $E_p(G)$ for $p\neq
2.$

\section{Orthogonal polynomials and Fourier series}

We give two standard applications for approximation of the
Szeg\H{o} kernel here. Since polynomials $\{p_n\}_{n=0}^{\infty}$
form a complete orthonormal system, we can restate \eqref{1.9} in
the following form:
\begin{equation} \label{2.1}
 \left(\sum_{k=n+1}^{\infty} |p_k(\zeta)|^2\right)^{1/2} \le C_1(G,F)\
n^{-\frac{\lambda}{4-2\lambda}}, \quad n\in\N,
\end{equation}
by  \eqref{1.1}. It is well known that $p_n$ converge to zero on
compact subsets of $G$ as $n\to\infty$ (cf. Chapter XVI of
\cite{Sz}). Equation \eqref{2.1} immediately gives an estimate for
the rate of decay of $p_n$ inside $G$.

\begin{theorem} \label{thm2.1}
Let $\partial G$ be piecewise analytic, with the smallest exterior
angle $\lambda\pi,\ 0<\lambda<2,$ at the junction points of the
analytic arcs. If $\zeta\in F$, where $F\subset G$ is compact,
then
\begin{equation} \label{2.2}
|p_n(\zeta)| \le C_1(G,F)\ n^{-\frac{\lambda}{4-2\lambda}}, \quad
n\in\N.
\end{equation}
\end{theorem}

The second application is related to the rates of convergence of
the Fourier series \eqref{1.3} for $f\in E_2(G)$ on compact
subsets of $G$. Observe that
\[
\|f\|_2 = \left(\sum_{k=0}^{\infty} |a_k|^2 \right)^{1/2}.
\]
Hence we have from \eqref{1.3} and Cauchy-Schwarz inequality that
\begin{align*}
\left|f(\zeta)-\sum_{k=0}^n a_k\, p_k(\zeta)\right| &=
\left|\sum_{k=n+1}^{\infty} a_k\, p_k(\zeta)\right| \le
\left(\sum_{k=n+1}^{\infty} |a_k|^2 \right)^{1/2}
\left(\sum_{k=n+1}^{\infty} |p_k(\zeta)|^2 \right)^{1/2} \\ &=
\left\|f-\sum_{k=0}^n a_k\, p_k\right\|_2
\left(\sum_{k=n+1}^{\infty} |p_k(\zeta)|^2 \right)^{1/2} \\ &\le
\|f\|_2 \left(\sum_{k=n+1}^{\infty} |p_k(\zeta)|^2 \right)^{1/2}.
\end{align*}
Thus we obtain the following result from \eqref{2.1}:

\begin{theorem} \label{thm2.2}
Let $\partial G$ be piecewise analytic, with the smallest exterior
angle $\lambda\pi,\ 0<\lambda<2,$ at the junction points of the
analytic arcs. Suppose that $f\in E_2(G)$ has the Fourier
expansion \eqref{1.3}. If $\zeta\in F$, where $F\subset G$ is
compact, then
\begin{equation} \label{2.3}
\left|f(\zeta)-\sum_{k=0}^n a_k\, p_k(\zeta)\right| \le C_1(G,F)\
n^{-\frac{\lambda}{4-2\lambda}}\ \left\|f-\sum_{k=0}^n a_k\,
p_k\right\|_2, \quad n\in\N.
\end{equation}
\end{theorem}

Results of this kind for domains with smooth boundaries were
previously proved by Szeg\H{o} \cite{Sz},  Rosenbloom and
Warschawski \cite{RW}, Smirnov and Lebedev \cite{SL}, Suetin
\cite{Su}, and others.

\section{Proofs}

{\bf Proof of Theorem \ref{thm1.1}.} Let $\psi:=\varphi^{-1}$. We
have
\begin{align*}
|\varphi(z)-J_{2n+1}(z)| &= \left|\int_{\zeta}^z
\left(\varphi'(t)-J_{2n+1}'(t)\right) dt\right| \\
&= \left|\int_{0}^{\varphi(z)} \left(\varphi'(\psi(u)) -
J_{2n+1}'(\psi(u))\right)\psi'(u)du\right| \\ &\le
\int_{0}^{\varphi(z)} \left|\left(\varphi'(\psi(u)) -
J_{2n+1}'(\psi(u))\right)\psi'(u)\right| |du|,
\end{align*}
where the integration is carried over the segment connecting $0$
and $\varphi(z)$ in $D$. Note that the function
$\left(\varphi'(\psi(u)) - J_{2n+1}'(\psi(u))\right)\psi'(u) = 1 -
J_{2n+1}'(\psi(u))\psi'(u)$ belongs to the Hardy class $H^1(D)$,
because $L$ is rectifiable. Hence we obtain by Fej\'er-Riesz
inequality (cf. Theorem 3.13 of \cite{Du}) that
\begin{align*}
|\varphi(z)-J_{2n+1}(z)| &\le \frac{1}{2} \int_{|u|=1}
\left|\left(\varphi'(\psi(u)) -
J_{2n+1}'(\psi(u)\right)\psi'(u)\right| |du| \\ &= \frac{1}{2}
\int_L \left|\varphi'(t) - J_{2n+1}'(t)\right| |dt|.
\end{align*}
Denote
\[
Q_n(z):=\left(\frac{l}{2\pi} \dis\sum_{k=0}^n |p_k(\zeta)|^2
\right)^{-1/2} \sum_{k=0}^n \overline{p_k(\zeta)} p_k(z),
\]
so that $J_{2n+1}'(z)=Q_n^2(z).$ We continue with this notation,
by using Cauchy-Schwarz and Minkowski inequalities:
\begin{align*}
|\varphi(z)-J_{2n+1}(z)| &\le \frac{1}{2} \int_L
\left|\left(\sqrt{\varphi'(t)}\right)^2 - Q_n^2(t)\right| |dt| \\
&= \frac{1}{2} \int_L \left|\sqrt{\varphi'(t)}-Q_n(t)\right|\,
\left|\sqrt{\varphi'(t)}+Q_n(t)\right| |dt|  \\ &\le \frac{l}{2}\
\|\sqrt{\varphi'}-Q_n\|_2 \ \|\sqrt{\varphi'}+Q_n\|_2 \\ &\le
\frac{l}{2}\ \|\sqrt{\varphi'}-Q_n\|_2 \left(\|\sqrt{\varphi'}\|_2
+ \|Q_n\|_2\right).
\end{align*}
Observe that
\[
\|\sqrt{\varphi'}\|_2 = \left(\frac{1}{l} \int_L |\varphi'(z)|\,
|dz| \right)^{1/2} = \left(\frac{1}{l} \int_{|w|=1} |dw|
\right)^{1/2} = \sqrt{\frac{2\pi}{l}},
\]
and that
\[
\|Q_n\|_2=\left(\frac{l}{2\pi} \dis\sum_{k=0}^n |p_k(\zeta)|^2
\right)^{-1/2} \left(\dis\sum_{k=0}^n |p_k(\zeta)|^2 \right)^{1/2}
= \sqrt{\frac{2\pi}{l}},
\]
by orthonormality of the polynomials $p_k$. Thus
\begin{align} \label{3.1}
\|\varphi-J_{2n+1}\|_{\infty} \le \sqrt{2\pi l}\
\|\sqrt{\varphi'}-Q_n\|_2.
\end{align}
We now estimate the norm on the right of the above inequality.
Recall that
\[
\sqrt{\varphi'(z)} = \sqrt{\frac{2\pi}{l}}\,
\frac{K(z,\zeta)}{\sqrt{K(\zeta,\zeta)}},
\]
by \eqref{1.6}. Therefore,
\begin{align*}
\sqrt{\varphi'(z)} - Q_n(z) &= \sqrt{\frac{2\pi}{l}} \left(
\frac{K(z,\zeta)}{\sqrt{K(\zeta,\zeta)}} - \left(\dis\sum_{k=0}^n
|p_k(\zeta)|^2 \right)^{-1/2} \sum_{k=0}^n \overline{p_k(\zeta)}
p_k(z) \right) \\
&= \sqrt{\frac{2\pi}{l K(\zeta,\zeta)}} \left(K(z,\zeta) -
\sum_{k=0}^n \overline{p_k(\zeta)} p_k(z) \right) \\ &+
\sqrt{\frac{2\pi}{l K(\zeta,\zeta)}} \left(1 -
\left(\frac{K(\zeta,\zeta)}{\sum_{k=0}^n |p_k(\zeta)|^2}
\right)^{1/2} \right) \sum_{k=0}^n \overline{p_k(\zeta)} p_k(z) .
\end{align*}
It follows that
\begin{align*}
\|\sqrt{\varphi'}-Q_n\|_2 &\le \sqrt{\frac{2\pi}{l
K(\zeta,\zeta)}} \left\|K(\cdot,\zeta) - \sum_{k=0}^n
\overline{p_k(\zeta)} p_k(\cdot) \right\|_2 \\ &+
\sqrt{\frac{2\pi}{l K(\zeta,\zeta)}}
\left(\left(\frac{K(\zeta,\zeta)}{\sum_{k=0}^n |p_k(\zeta)|^2}
\right)^{1/2} - 1\right) \left(\sum_{k=0}^n |p_k(\zeta)|^2
\right)^{1/2} \\ &= \sqrt{\frac{2\pi}{l K(\zeta,\zeta)}}
\left\|K(\cdot,\zeta) - \sum_{k=0}^n \overline{p_k(\zeta)}
p_k(\cdot) \right\|_2 \\ &+ \sqrt{\frac{2\pi}{l K(\zeta,\zeta)}}
\left(\left(K(\zeta,\zeta)\right)^{1/2} - \left(\sum_{k=0}^n
|p_k(\zeta)|^2 \right)^{1/2} \right).
\end{align*}
Since
\begin{align*}
\left(K(\zeta,\zeta)\right)^{1/2} - \left(\sum_{k=0}^n
|p_k(\zeta)|^2 \right)^{1/2} &= \left\|K(\cdot,\zeta)\right\|_2 -
\left\|\sum_{k=0}^n \overline{p_k(\zeta)} p_k(\cdot) \right\|_2
\\ &\le \left\|K(\cdot,\zeta) - \sum_{k=0}^n \overline{p_k(\zeta)}
p_k(\cdot) \right\|_2,
\end{align*}
we obtain that
\begin{align*}
\|\sqrt{\varphi'}-Q_n\|_2 &\le 2 \sqrt{\frac{2\pi}{l
K(\zeta,\zeta)}} \left\|K(\cdot,\zeta) - \sum_{k=0}^n
\overline{p_k(\zeta)} p_k(\cdot) \right\|_2.
\end{align*}
Equation \eqref{3.1} now gives that
\begin{align} \label{3.2}
\|\varphi-J_{2n+1}\|_{\infty} \le
\frac{4\pi}{\sqrt{K(\zeta,\zeta)}} \left\|K(\cdot,\zeta) -
\sum_{k=0}^n \overline{p_k(\zeta)} p_k(\cdot) \right\|_2.
\end{align}
Recall that
\[
K(\zeta,\zeta)= \frac{l}{2\pi}\, \varphi'(\zeta),
\]
by \eqref{1.5}. Using Corollary 1.4 of \cite[p. 9]{Po}, we have
\[
|\varphi'(\zeta)| \ge \frac{1}{4}\,
\frac{1-|\varphi(\zeta)|^2}{\textup{dist}(\zeta,L)} =
\frac{1}{4\,\textup{dist}(\zeta,L)},
\]
where $\textup{dist}(\zeta,L)$ is the distance from $\zeta$ to
$L.$ Clearly, $l \ge 2\pi\,\textup{dist}(\zeta,L)$, so that
\[
K(\zeta,\zeta)= \frac{l\,\varphi'(\zeta)}{2\pi} \ge
\varphi'(\zeta)\,\textup{dist}(\zeta,L) \ge \frac{1}{4}.
\]
Combining this with \eqref{3.2}, we obtain \eqref{1.8}.

\db

{\bf Proof of Theorem \ref{thm1.2}.} We start by recalling that
the partial sum $\sum_{k=0}^n \overline{p_k(\zeta)} p_k(z)$ for
$K(z,\zeta)$ is its best approximation in $E_2(G)$ (with $\zeta\in
G$ fixed) among all polynomials of degree at most $n.$ Thus we
construct a sequence of polynomials with good approximative
properties, which gives the desired upper bound \eqref{1.9}.

It is clear from \eqref{1.4} that approximation of $K(z,\zeta)$ is
equivalent to approximation of $\sqrt{\varphi'(z)}.$ We use a
method resembling that of Andrievskii and Gaier \cite{AG}. We
first continue the mapping $\varphi$ conformally beyond the
boundary $L$, by using reflections across the analytic arcs $L_i,\
L=\cup_{i=1}^m L_i.$ Suppose that $\tau_i$ is a mapping such that
$L_i=\tau_i([0,1])$, which is conformal in an open neighborhood of
$[0,1].$ Then we can find a symmetric lens shaped domain $S_i$,
bounded by two circular arcs subtended by $[0,1]$, whose closure
is contained in this open neighborhood of $[0,1].$ Defining
\[
\tilde G := G \cup \left( \cup_{i=1}^m \tau_i(S_i) \right),
\]
we extend $\varphi$ into $\tilde G$ as follows:
\[
\varphi (z) :=\frac{1}{\overline{\varphi \left[ \tau_i \left(
\overline{\tau^{-1}_i(z)} \right) \right]} }, \qquad z \in
\tau_i(S_i) \backslash \overline{G},
\]
where $i=1,\ldots,m.$ The boundary $\partial\tilde G$ consists of
$m$ analytic arcs $\Gamma_i$ that share endpoints with the arcs
$L_i$ of $\partial G$:
\[
\partial\tilde G \cap \partial G = \{z_i\}_{i=1}^m,
\]
which are clearly the corner points of $\partial G.$ Since each
$\tau_i,\ i=1,\ldots,m,$ is conformal and has bounded derivative
(together with its inverse) on $S_i$, we obtain the inequalities
\begin{align} \label{3.3}
{\rm dist}(z,\partial G) \ge c_1 \min_{1\le i\le m} |z-z_i|,
\qquad z\in \partial\tilde G,
\end{align}
where ${\rm dist}(z,\partial G)$ is the distance from $z$ to
$\partial G$, and
\begin{align} \label{3.4}
|\gamma| \le c_2 |z-t|, \qquad z,t\in
\partial\tilde G,
\end{align}
where $|\gamma|$ is the length of the shorter arc $\gamma\subset
\partial\tilde G,$ connecting $z$ and $t$. We denote various
positive constants by $c_1,c_2,$ etc.

Let $\Gamma_j$ be an arc of $\partial \tilde G$, with the
endpoints $z_j$ and $z_{j+1}$, and let $\zeta_j \in \Gamma_j$ be a
fixed point, $j=1,\ldots,m.$ Note that $\zeta_j$ divides
$\Gamma_j$ into $\Gamma_j^1$ and $\Gamma_j^2$, so that $\partial
\tilde G = \bigcup_{j =1}^m \bigcup_{i =1}^2 \Gamma_j^i$. We
obtain from Cauchy's integral formula for the continuation of
$\sqrt{\varphi'}$ into $\tilde G$ that
\begin{align} \label{3.5}
\sqrt{\varphi'(z)} = \frac{1}{2 \pi i}  \int_{\partial \tilde G}
\frac{\sqrt{\varphi'(t)}}{t-z}\, dt = \frac{1}{2 \pi i} \sum_{j
=1}^m \sum_{i=1}^2 \int_{\Gamma_j^i}
\frac{\sqrt{\varphi'(t)}}{t-z}\, dt, \quad z \in \tilde G.
\end{align}
Hence the problem is reduced to approximation of functions of the
form
\begin{align} \label{3.6}
g(z) := \int_{\gamma} \frac{\sqrt{\varphi'(t)}}{t-z}\, dt
\end{align}
in $E_2(G)$ norm, where $\gamma$ is any of the arcs $\Gamma_j^i,$
with $i=1,2$ and $j=1,\ldots,m.$

Let $\Omega:=\overline{\C}\setminus\overline{G}.$ Consider the
standard conformal mapping $\Phi:\Omega\to\Delta,$ where
$\Delta:=\{w:|w|>1\},$ normalized by $\Phi(\infty)=\infty$ and
$\Phi'(\infty)>0.$ We define the level curves of $\Phi$ by
\[
L_n:=\{z: |\Phi(z)|=1+1/n\}, \qquad n\in\N.
\]
Denote by $\gamma_2$ the part of $\gamma$ from its endpoint
$\zeta_j\in\Gamma_j$ to the first point $\xi$ of intersection with
$L_n$, so that $\gamma_2 \subset \{z: |\Phi(z)|>1+1/n\}.$ Then
$\gamma_1:=\gamma\setminus\gamma_2$ connects $\xi$ with the corner
point $z_j$ of $L.$ Write
\begin{align} \label{3.7}
g(z) := \int_{\gamma_1} \frac{\sqrt{\varphi'(t)}}{t-z}\, dt +
\int_{\gamma_2} \frac{\sqrt{\varphi'(t)}}{t-z}\, dt =: g_1(z) +
g_2(z).
\end{align}
We show that $\|g_1\|_2\to 0$ sufficiently fast as $n\to\infty,$
while $g_2$ is well approximated by polynomials of degree $n.$ To
estimate the norm of $g_1$, we need to know the behavior of
$\sqrt{\varphi'}$ near the corner point $z_j\in L.$ This is
conveniently found from the asymptotic expansion of Lehman
\cite{Le}. Assume that $z_j=0$ and that $\lambda_j\pi,\
0<\lambda_j<2,$ is the exterior angle formed by $L$ at this point.
Then we have in a neighborhood of $z_j=0$ that
\[
\varphi(z)-\varphi(0) = b\ z^{\frac{1}{2-\lambda_j}} +
o\left(z^{\frac{1}{2-\lambda_j}}\right) \qquad \mbox{as } z \to 0,
\]
where $b\neq 0$,
and
\[
\varphi'(z) = \frac{b}{2-\lambda_j}\ z^{\frac{1}{2-\lambda_j}-1} +
o\left(z^{\frac{1}{2-\lambda_j}-1}\right) \qquad \mbox{as } z \to
0.
\]
Hence there exists a constant $c_3>0$ such that
\begin{align} \label{3.8}
\left|\sqrt{\varphi'(z)}\right| \le c_3\ |z|^{\alpha}, \qquad z\in
\tilde G \cup \partial \tilde G,
\end{align}
where we set
\[
\alpha:=\frac{1}{4-2\lambda_j}-\frac{1}{2}.
\]
For the endpoints $\xi\in L_n$ and $0$ of $\gamma_1$, we let
\[
d_n:=|\xi-0|=|\xi|.
\]
It follows from \eqref{3.4} that
\[
|\gamma_1| \le c_2 d_n.
\]
We now estimate that
\begin{align} \label{3.9}
\|g_1\|_2^2 = \frac{1}{l} \int_L \left| \int_{\gamma_1}
\frac{\sqrt{\varphi'(t)}}{t-z}\, dt \right|^2\ |dz| \le c_4 \int_L
\left( \int_{\gamma_1} \frac{|t|^{\alpha} |dt|}{|t-z|} \right)^2\
|dz|,
\end{align}
by \eqref{3.7} and \eqref{3.8}. Note that if $z\in L$ satisfies
$|z|\ge d_n$, then $|t-z|\sim |z|$ by \eqref{3.3}. Consequently,
\begin{align} \label{3.10}
\int_{L\cap\{|z|\ge d_n\}} \left( \int_{\gamma_1}
\frac{|t|^{\alpha} |dt|}{|t-z|} \right)^2\ |dz| \le c_5
\int_{L\cap\{|z|\ge d_n\}} \left( \frac{d_n^{\alpha+1}}{|z|}
\right)^2\ |dz| \le c_6\,d_n^{2\alpha+1}.
\end{align}
On the other hand, if $z\in L$ satisfies $|z|\le d_n$, then
$|t-z|\sim |t|+|z|$ by \eqref{3.3}, and we obtain with help of
\eqref{3.4} that
\begin{align} \label{3.11}
\int_{L\cap\{|z|\le d_n\}} \left( \int_{\gamma_1}
\frac{|t|^{\alpha} |dt|}{|t-z|} \right)^2\ |dz| &\le c_7
\int_0^{c_8 d_n} \left( \int_0^{c_9 d_n} \frac{s^{\alpha}ds}{s+r}
\right)^2\ dr \\ \nonumber &\le c_7 \int_0^{c_8 d_n} \left(
\int_0^r \frac{s^{\alpha}}{r}\,ds + \int_r^{c_9 d_n}
s^{\alpha-1}ds \right)^2\ dr \\ \nonumber &= c_7 \int_0^{c_8 d_n}
\left( \frac{r^{\alpha}}{\alpha+1} +
\frac{(c_9d_n)^{\alpha}-r^{\alpha}}{\alpha} \right)^2\ dr \\
\nonumber &\le c_{10}\, d_n^{2\alpha+1},
\end{align}
for $\alpha\neq 0$. If $\alpha=0$ then we estimate
\begin{align*}
\int_{L\cap\{|z|\le d_n\}} \left( \int_{\gamma_1}
\frac{|dt|}{|t-z|} \right)^2\ |dz| &= \int_{L\cap\{|z|\le d_n\}}
\left( \int_{\gamma_1} \frac{|t|^{1/2}\, |t|^{-1/2}\, |dt|}{|t-z|}
\right)^2\ |dz| \\&\le c_2\,d_n \int_{L\cap\{|z|\le d_n\}} \left(
\int_{\gamma_1} \frac{|t|^{-1/2}\, |dt|}{|t-z|} \right)^2\ |dz|
\\&\le c_2\,d_n\, c_{10}\, d_n^{2(-1/2)+1}= c_2\,c_{10}\, d_n,
\end{align*}
as above. Combining \eqref{3.9}-\eqref{3.11}, we have that
\begin{align} \label{3.12}
\|g_1\|_2 \le c_{11}\, d_n^{\alpha+1/2} \le c_{11}\,
d_n^{\frac{1}{4-2\lambda}},
\end{align}
where $\lambda=\min_{1\le j \le m} \lambda_j.$

The next step is the construction of approximating polynomials
$P_n$ for $g_2$. This is accomplished by using Dzjadyk kernels
(see, e.g., \cite{ABD}) of the form
\[
K_n(t,z)=\sum_{i=0}^n a_i(t) z^i, \qquad n\in\N,
\]
which approximate the Cauchy kernel. It was proved in Lemma 5 of
\cite{AG} that a sequence of such kernels can be selected, so that
for any fixed $m\in\N$, and for all $t\in\gamma$ with
$|\Phi(t)|\ge 1+1/n,$ we have
\begin{align} \label{3.13}
\left|\frac{1}{t-z} - K_n(t,z)\right| \le c_{12}\,
\frac{d_n^m}{|t-z|^{m+1}}, \qquad z\in L,
\end{align}
for all sufficiently large $n\in\N.$ In particular, \eqref{3.13}
holds for $t\in\gamma_2.$ Define the polynomials
\[
P_n(z):=\int_{\gamma_2} \sqrt{\varphi'(t)}\, K_n(t,z)\,dt,
\]
and estimate
\begin{align*}
\|g_2-P_n\|_2^2 &= \frac{1}{l} \int_L \left| \int_{\gamma_2}
\left( \frac{1}{t-z}-K_n(t,z) \right) \sqrt{\varphi'(t)}\, dt
\right|^2\ |dz| \\ &\le c_{13} d_n^{2m} \int_L \left(
\int_{\gamma_2} \frac{|t|^{\alpha} |dt|}{|t-z|^{m+1}} \right)^2\
|dz|,
\end{align*}
by \eqref{3.13} and \eqref{3.8}. Observe that $|t-z| \sim |t|+|z|$
for $t\in\gamma_2.$ Therefore, we have for $m>\alpha$ that
\begin{align*}
\int_L \left( \int_{\gamma_2} \frac{|t|^{\alpha}
|dt|}{|t-z|^{m+1}} \right)^2\ |dz| &\le c_{14} \int_0^{c_{15}}
\left( \int_{c_{16} d_n}^{c_{17}} \frac{s^{\alpha}ds}{(s+r)^{m+1}} \right)^2\ dr \\
&\le c_{14} \int_0^{c_{16} d_n} \left( \int_{c_{16} d_n}^{c_{17}}
s^{\alpha-m-1}\,ds \right)^2\ dr \\ &+ c_{14} \int_{c_{16}
d_n}^{c_{15}} \left(r^{-m-1} \int_{c_{16} d_n}^r
s^{\alpha}\,ds + \int_r^{c_{17}} s^{\alpha-m-1}\,ds \right)^2\ dr \\
&\le c_{18}\, d_n^{2(\alpha-m)+1} + c_{19}
\int_{c_{16}d_n}^{c_{15}} r^{2(\alpha-m)}\,dr \\
&\le c_{20}\, d_n^{2(\alpha-m)+1}.
\end{align*}
It follows that
\begin{align} \label{3.14}
\|g_2-P_n\|_2 \le c_{21}\, d_n^{\alpha+1/2} \le c_{21}\,
d_n^{\frac{1}{4-2\lambda}}.
\end{align}
Collecting \eqref{3.12}, \eqref{3.14} and \eqref{3.7} together, we
obtain
\begin{align} \label{3.15}
\|g-P_n\|_2 \le \|g_1\|_2 + \|g_2-P_n\|_2 \le c_{22}\,
d_n^{\frac{1}{4-2\lambda}}.
\end{align}

Recall that $d_n=|\xi|,$ where $\xi\in L_n \cap \gamma_1.$
Applying the results of \cite{Le} to the conformal mapping
$\Psi:=\Phi^{-1},$ we obtain
\[
z=\Psi(\Phi(z))-\Psi(\Phi(0))=a
\left(\Phi(z)-\Phi(0)\right)^{\lambda_j} +
o\left(\left(\Phi(z)-\Phi(0)\right)^{\lambda_j}\right)
\quad\mbox{as }z\to 0,
\]
where $\lambda_j\pi$ is the exterior angle at $z_j=0,$ and $a\neq
0.$ Thus
\[
d_n=|\xi| \le c_{23} \min_{z\in L_n} |z| \le c_{24}\,
n^{-\lambda_j} \le c_{24}\, n^{-\lambda}, \qquad n\in\N,
\]
and
\begin{align} \label{3.16}
\|g-P_n\|_2 \le c_{25}\, n^{-\frac{\lambda}{4-2\lambda}}, \qquad
n\in\N,
\end{align}
by \eqref{3.15}. Hence there exists a sequence of polynomials
$Q_n$ such that
\begin{align} \label{3.17}
\|\sqrt{\varphi'}-Q_n\|_2 \le c_{26}\,
n^{-\frac{\lambda}{4-2\lambda}}, \qquad n\in\N.
\end{align}
Since
\[
K(z,\zeta)=\frac{l}{2\pi} \sqrt{\varphi'(z)\varphi'(\zeta)},
\]
we obtain \eqref{1.9} from the previous estimate:
\[
\left\|K(\cdot,\zeta)-\frac{l}{2\pi} \sqrt{\varphi'(\zeta)}\ Q_n
\right\|_2 \le \frac{c_{26}l}{2\pi} \sqrt{\varphi'(\zeta)}\
n^{-\frac{\lambda}{4-2\lambda}}, \qquad n\in\N,
\]
where $\zeta\in G$ was fixed throughout this proof. We now lift
this restriction and allow $\zeta$ vary within a compact set
$F\subset G.$ Note that $|\varphi'(\zeta)|$ is uniformly bounded
on $F$, and so are other constants in the above proof. Indeed, we
only need to verify this for the constant $c_3$ of \eqref{3.8},
arising from the first term in the expansion of $\varphi'(z).$ One
can obtain a conformal mapping $\tilde\varphi$ of $G$ onto $D$
with $\tilde\varphi(\tilde\zeta)=0$ and
$\tilde\varphi'(\tilde\zeta)>0$, for any $\tilde\zeta\in F$, by
composing $\varphi$ with a M\"obius self-map of the unit disk,
which is conformal in an open neighborhood of $\overline D.$ It
follows that \eqref{3.8} holds for all such mappings
$\tilde\varphi$ with a constant $c_3$ uniformly bounded for
$\tilde\zeta\in F$.

\db

{\bf Proof of Theorem \ref{thm1.4}.} We proceed as in the proof of
Theorem \ref{thm1.1}, denoting
\[
Q_n(z):=\left(\frac{l}{2\pi} \dis\sum_{k=0}^n |p_k(\zeta)|^2
\right)^{-1/2} \sum_{k=0}^n \overline{p_k(\zeta)} p_k(z),
\]
so that $J_{2n+1}'(z)=Q_n^2(z).$ It follows from \eqref{1.6} that
\begin{align*}
\sqrt{\varphi'(z)} - Q_n(z) &= \sqrt{\frac{2\pi}{l}} \left(
\frac{K(z,\zeta)}{\sqrt{K(\zeta,\zeta)}} - \left(\dis\sum_{k=0}^n
|p_k(\zeta)|^2 \right)^{-1/2} \sum_{k=0}^n \overline{p_k(\zeta)}
p_k(z) \right) \\
&= \sqrt{\frac{2\pi}{l K(\zeta,\zeta)}} \left(K(z,\zeta) -
\sum_{k=0}^n \overline{p_k(\zeta)} p_k(z) \right) \\ &+
\sqrt{\frac{2\pi}{l K(\zeta,\zeta)}} \frac{\left(\sum_{k=0}^n
|p_k(\zeta)|^2 \right)^{1/2} -
\left(K(\zeta,\zeta)\right)^{1/2}}{\left(\sum_{k=0}^n
|p_k(\zeta)|^2 \right)^{1/2}} \sum_{k=0}^n \overline{p_k(\zeta)}
p_k(z) .
\end{align*}
Recall that
\[
\lim_{n\to\infty} \sum_{k=0}^n \overline{p_k(\zeta)} p_k(z) =
K(z,\zeta),
\]
where convergence is uniform for $z,\zeta\in F,$ and that
\[
0<c_1<|K(z,\zeta)|<c_2<+\infty, \qquad z,\zeta\in F.
\]
Hence
\begin{align*}
\max_{z\in F}|\sqrt{\varphi'(z)} - Q_n(z)| &\le c_3 \max_{z\in F}
\left|K(z,\zeta) - \sum_{k=0}^n \overline{p_k(\zeta)} p_k(z)
\right| \\ &+ c_4 \left(\left(K(\zeta,\zeta)\right)^{1/2} -
\left(\sum_{k=0}^n |p_k(\zeta)|^2 \right)^{1/2}\right).
\end{align*}
The first term is estimated by \eqref{1.9}-\eqref{2.1} and
Cauchy-Schwarz inequality:
\begin{align*}
\max_{z,\zeta\in F} \left|K(z,\zeta) - \sum_{k=0}^n
\overline{p_k(\zeta)} p_k(z) \right| &= \max_{z,\zeta\in F}
\left|\sum_{k=n+1}^{\infty} \overline{p_k(\zeta)} p_k(z) \right|\\
&\le \left(\sum_{k=n+1}^{\infty} |p_k(z)|^2 \right)^{1/2}
\left(\sum_{k=n+1}^{\infty} |p_k(\zeta)|^2 \right)^{1/2}\\ &\le
(C_1(G,F))^2\ n^{-\frac{\lambda}{2-\lambda}}, \quad n\in\N.
\end{align*}
Thus for the second term we also have
\begin{align*}
\left(K(\zeta,\zeta)\right)^{1/2} - \left(\sum_{k=0}^n
|p_k(\zeta)|^2 \right)^{1/2} &\le c_5 \left(K(\zeta,\zeta) -
\sum_{k=0}^n |p_k(\zeta)|^2 \right) \\ &\le c_5\,(C_1(G,F))^2\
n^{-\frac{\lambda}{2-\lambda}}, \quad n\in\N.
\end{align*}
Combining these estimates, we obtain that
\begin{align*}
\max_{z\in F}|\sqrt{\varphi'(z)} - Q_n(z)| \le c_6\,
n^{-\frac{\lambda}{2-\lambda}}, \qquad n\in\N.
\end{align*}
It immediately follows that
\begin{align*}
\max_{z\in F}|\varphi'(z) - Q_n^2(z)| &\le \max_{z\in
F}|\sqrt{\varphi'(z)} + Q_n(z)|\ \max_{z\in F}|\sqrt{\varphi'(z)}
- Q_n(z)| \\ &\le c_7\, n^{-\frac{\lambda}{2-\lambda}}, \qquad
n\in\N,
\end{align*}
and that
\begin{align*}
\max_{z\in F}|\varphi(z) - J_{2n+1}(z)| \le \max_{z\in F}
\int_{\zeta}^z |\varphi'(t) - Q_n^2(t)|\,|dt| \le c_8\,
n^{-\frac{\lambda}{2-\lambda}}, \quad n\in\N,
\end{align*}
where all constants in this proof are independent of $z,\zeta\in
F$ and $n\in\N.$

\db

{\it Igor E. Pritsker}  \hskip 2in {\sc E-mail:}
igor@math.okstate.edu

{\sc Address:} Department of Mathematics, 401 Mathematical
Sciences, Oklahoma State University, Stillwater, OK 74078-1058,
U.S.A.

\end{document}